\documentclass[aos,nameyear,dvips]{arximspdf}
\usepackage{graphicx}
\doi{10.1214/09-AOS786}
\volume{38}
\issue{6}
\pubyear{2010}
\firstpage{3458}
\lastpage{3486}

\makeatletter

\newcommand{\cov}{\operatorname{cov}}
\newcommand{\var}{\operatorname{var}}
\newcommand{\mt}{\mathcal{T}}

\newtheorem{theorem}{Theorem}
\newtheorem{lem}{Lemma}

\newcommand{\T}{\mathcal{T}}

\makeatother

\begin{document}
\begin{frontmatter}

\title{Empirical dynamics for longitudinal data}
\runtitle{Empirical dynamics}

\begin{aug}
\author[A]{\fnms{Hans-Georg} \snm{M\"uller}\corref{}\thanksref{t1}\ead[label=e1]{mueller@wald.ucdavis.edu}} and
\author[B]{\fnms{Fang} \snm{Yao}\thanksref{t2}\ead[label=e2]{fyao@utstat.toronto.edu}}
\runauthor{H.-G. M\"uller and F. Yao}
\affiliation{University of California, Davis and University of Toronto}
\address[A]{Department of Statistics \\
University of California, Davis \\
One Shields Avenue \\
Davis, California 95616\\
USA \\
\printead{e1}} %adresu isvedimo komanda gale!
\address[B]{Department of Statistics\\
University of Toronto\\
100 St. George Street\\
Toronto, ON M5S 3G3\\
Canada\\
\printead{e2}}
\end{aug}

\thankstext{t1}{Supported in part by NSF Grant DMS-08-06199.}
\thankstext{t2}{Supported in part by a NSERC Discovery grant.}

% HISTORY:
\received{\smonth{7} \syear{2009}}
\revised{\smonth{11} \syear{2009}}

% ABSTRACT
%
\begin{abstract}
We demonstrate that the processes underlying on-line auction
price bids and many other longitudinal data can be represented by an
empirical first order stochastic ordinary differential equation with
time-varying coefficients and a smooth drift process. This equation
may be empirically obtained from longitudinal observations for a
sample of subjects and does not presuppose specific knowledge of the
underlying processes. For the nonparametric estimation of the
components of the differential equation, it suffices to have
available sparsely observed longitudinal measurements which may be
noisy and are generated by underlying smooth random trajectories
for each subject or experimental unit in the sample. The drift
process that drives the equation determines how closely individual
process trajectories follow a deterministic approximation of the
differential equation. We provide estimates for trajectories and
especially the variance function of the drift process. At each
fixed time point, the proposed empirical dynamic model implies a
decomposition of the derivative of the process underlying the
longitudinal data into a component explained by a linear component
determined by a varying coefficient function dynamic equation and an
orthogonal complement that corresponds to the drift process. An
enhanced perturbation result enables us to obtain improved
asymptotic convergence rates for eigenfunction derivative estimation
and consistency for the varying coefficient function and the
components of the drift process. We illustrate the differential
equation with an application to the dynamics of on-line auction
data.\looseness=-1
\end{abstract}

% KEYWORDS
%

\end{frontmatter}

%s1 ###
\section{Introduction}
Recently, there has been increasing interest in analyzing on-line
auction data and in inferring the underlying dynamics that drive the
bidding process. Each series of price bids for a given auction
corresponds to pairs of random bidding times and corresponding bid
prices generated whenever a bidder places a bid
[Jank and Shmueli (\citeyear{jank051}, \citeyear{jank062}), \citet{jank08},
\citet{redd06}]. Related longitudinal
data where similar sparsely and irregularly sampled noisy
measurements are obtained are abundant in the social and life
sciences; for example, they arise in longitudinal growth studies.
While more traditional approaches of functional data analysis
require fully or at least densely observed trajectories [\citet{kirk89},
\citet{rams05}, \citet{gerv05}], more recent extensions
cover the case of sparsely
observed and noise-contaminated longitudinal data [\citet{mull054},
\citet{wang052}].

A common assumption of approaches for longitudinal data grounded in
functional data analysis is that such data are generated by an
underlying smooth and square integrable stochastic process
[\citet{sy97}, \citet{stan98}, \citet{rice04},
\citet{zhao04}, \citet{mull067}].
The derivatives of the trajectories of such
processes are central for assessing the dynamics of the underlying
processes [\citet{rams00}, \citet{mas07}]. Although this
is difficult for
sparsely recorded data, various approaches for estimating
derivatives of individual trajectories nonparametrically by pooling
data from samples of curves and using these derivatives for
quantifying the underlying dynamics have been developed
[\citet{mull842}, \citet{reit08}, \citet{wang081},
\citet{wang082}]. Related work on
nonparametric methods for derivative estimation can be found in
\citet{mull841}, \citet{hard85} and on the role of
derivatives for the
functional linear model in \citet{mas09}.

We expand here on some of these approaches and investigate an
empirical dynamic equation. This equation is distinguished from
previous models that involve differential equations in that it is
empirically determined from a sample of trajectories, and does not
presuppose knowledge of a specific parametric form of a
differential equation which generates the data, except that we
choose it to be a first order equation. This stands in contrast to
current approaches of modeling dynamic systems, which are
``parametric'' in the sense that a prespecified differential
equation is assumed. A typical example for such an approach has
been developed by \citet{rams07}, where a prior specification of a
differential equation is used to guide the modeling of the data,
which is done primarily for just one observed trajectory. A
problem with parametric approaches is that diagnostic tools to
determine whether these equations fit the data either do not
exist, or where they do, are not widely used, especially as
nonparametric alternatives to derive differential equations have
not been available. This applies especially to the case where one
has data on many time courses available, providing strong
motivation to explore nonparametric approaches to quantify
dynamics. Our starting point is a nonparametric approach to
derivative estimation by local polynomial fitting of the
derivative of the mean function and of partial derivatives of the
covariance function of the process by pooling data across all
subjects [\citet{mull091}].

We show that each trajectory satisfies a first order stochastic
differential equation where the random part of the equation
resides in an additive smooth drift process which drives the
equation; the size of the variance of this process determines to
what extent the time evolution of a specific trajectory is
determined by the nonrandom part of the equation over various
time subdomains, and therefore is of tantamount interest. We
quantify the size of the drift process by its variance as a
function of time. Whenever the variance of the drift process $Z$
is small relative to the variance of the process $X$, a
deterministic version of the differential equation is particularly
useful as it then explains a large fraction of the variance of the
process.

The empirical stochastic differential equation can be easily
obtained for various types of longitudinal data. This approach thus
provides a novel perspective to assess the dynamics of longitudinal
data and permits insights about the underlying forces that shape the
processes generating the observations, which would be hard to obtain
with other methods. We illustrate these empirical dynamics by
constructing the stochastic differential equations that govern
online auctions with sporadic bidding patterns.

We now describe a data model for longitudinally collected
observations, which reflects that the data consist of sparse,
irregular and noise-corrupted measurements of an underlying smooth
random trajectory for each subject or experimental unit
[\citet{mull054}], the dynamics of which is of interest. Given $n$
realizations $X_i$ of the underlying process $X$ on a domain $\T$
and $N_i$ of an integer-valued bounded random variable $N$, we
assume that $N_i$ measurements $Y_{ij}$, $i = 1,\ldots, n, j =
1,\ldots,N_i$, are obtained at random times $T_{ij}$, according to
%
%e1 ###
%
\begin{equation} \label{kl2}
Y_{ij} = %Y_i(T_{ij}) =
X_i(T_{ij})
+ \varepsilon_{ij},\qquad
%=\mu(T_{ij}) +
%+ \varepsilon_{ij},
T_{ij} \in\T,
\end{equation}
where $\varepsilon_{ij}$ are zero mean i.i.d. measurement errors,
with $\var(\varepsilon_{ij}) = \sigma^2$, independent of all other
random components.

The paper is organized as follows. In Section~\ref{sec2}, we review
expansions in eigenfunctions and functional principal components, which
we use directly as the basic tool for dimension
reduction---alternative implementations with B-splines or P-splines could also be
considered [\citet{shi96}, \citet{rice01}, \citet
{yao06}]. We also introduce the empirical stochastic differential
equation and discuss the decomposition of variance it entails.
Asymptotic properties of estimates for the components of the
differential equation, including variance function of the drift
process, coefficient of determination associated with the dynamic
system and auxiliary results on improved rates of convergence for
eigenfunction derivatives are the theme of Section~\ref{sec3}.
Background on related perturbation results can be found in
\citet{daux82}, \citet{fine87}, \citet{kato95},
\citet{mas03}. Section~\ref{sec4} contains the illustration of the
differential equation with auction data, followed by a brief discussion
of some salient features of the proposed approach in Section
\ref{sec5}. Additional discussion of some preliminary formulas is
provided in Appendix~\ref{secA1}, estimation procedures are described
in Appendix~\ref{secA2}, assumptions and auxiliary results are in
Appendix~\ref{secA3} and proofs in Appendix~\ref{secA4}.

%s2 ###
\section{Empirical dynamics}\label{sec2}

%s2.1 ###
\subsection{Functional principal components and eigenfunction derivatives}

$\!\!$A~key methodology for dimension reduction and modeling of the
underlying stochastic processes $X$ that generate the longitudinal
data, which usually are sparse, irregular and noisy as in
(\ref{kl2}), is Functional Principal Component Analysis (FPCA).
Processes are assumed to be square integrable with mean function
$E\{X(t)\} = \mu(t)$ and auto-covariance function $\cov\{X(t),X(s)\}
= G(t, s)$, $s,t \in\mt$, which is smooth, symmetric and
nonnegative definite. Using $G$ as kernel in a linear operator
leads to the Hilbert--Schmidt operator $(\mathbf{G}f)(t) = \int_\mt G(t,
s)f(s)\,ds$. We denote the ordered eigenvalues (in declining order) of
this operator by $\lambda_1 > \lambda_2 > \cdots\geq0$ and the
corresponding orthonormal eigenfunctions by $\phi_k(t)$. We assume
that all eigenvalues are of multiplicity $1$ in the sequel.

It is well known that the kernel $G$ has the representation $G(t, s)
=\break \sum_{k=1}^{\infty}\lambda_k\times\phi_k(t)\phi_k(s)$ and the
trajectories generated by the process satisfy the
Karhunen--Lo\`{e}ve representation
[\citet{gren50}] $X_i(t) = \mu(t) +
\sum_{k=1}^{\infty}\xi_{ik}\phi_{k}(t)$. Here the $\xi_{ik} =
\int_{\mt}\{X_i(t)-\mu(t)\}\phi_{k}(t)\,dt$, $k=1,2,\ldots,$
$i=1,\ldots,\break n$, are the functional principal components (FPCs) of
the random trajectories $X_i$. The $\xi_k$ are uncorrelated random
variables with $E(\xi_{k}) = 0$ and $E\xi^2_{k} = \lambda_{k}$, with
$\sum_k \lambda_k < \infty$. Upon differentiating both sides, one
obtains
%
%e2 ###
%
\begin{equation}\label{kl3}
X_i^{(1)}(t) = \mu^{(1)}(t) + \sum_{k=1}^{\infty}\xi_{ik}\phi_k^{(1)}(t),
\end{equation}
where $\mu^{(1)}(t)$ and $\phi_k^{(1)}(t)$ are the derivatives of
mean and eigenfunctions.

The eigenfunctions $\phi_k$ are the solutions of the eigen-equations
$\int G(t,s)\times \phi_k(s)\, ds=\lambda_k\phi_k(t)$, under the constraint
of orthonormality. Under suitable regularity conditions, one
observes
%
%e3 ###
%
\begin{eqnarray}\label{kl4}
\frac{d}{dt}\int_{\mathcal{T}}G(t, s)
\phi_k(s)\,ds &=& \lambda_k\,\frac{d}{dt}\,\phi_k(t),\nonumber\\[-8pt]\\[-8pt]
\phi_k^{(1)}(t) &=&
\frac{1}{\lambda_k}\int_{\mt}
\frac{\partial}{\partial t}\,G(t, s)\phi_k(s)\,ds,\nonumber
\end{eqnarray}
which motivates corresponding eigenfunction
derivative estimates. A useful representation is
%
%e4 ###
%
\begin{eqnarray}\label{rep}
\cov\bigl\{X^{(\nu_1)}(t), X^{(\nu_2)}(s)\bigr\}=\sum_{k=1}^\infty
\lambda_k\phi_k^{(\nu_1)}(t) \phi_k^{(\nu_2)}(s),\nonumber\\[-8pt]\\[-8pt]
\eqntext{\nu_1, \nu_2
\in\{0,1\}, s,t \in\mt,}
\end{eqnarray}
which is an immediate consequence of the basic properties
of the
functional principal components $\xi_k$. For more details and
discussion, we refer to Appendix~\ref{secA1}.

It is worthwhile to note that the representation (\ref{kl3}) does
not correspond to
the Karhunen--Lo\`eve representation of the derivatives, which would
be based on orthonormal eigenfunctions of a linear Hilbert--Schmidt
operator defined
by the covariance kernel $\cov\{X^{(1)}(t), X^{(1)}(s)\}$. A method
to obtain this representation might proceed by first estimating $\cov\{
X^{(1)}(t), X^{(1)}(s)\}$
using (\ref{rep}) for
$\nu_1=\nu_2=1$ and suitable estimates $\phi_k^{(1)}$ for
eigenfunction derivatives, then directly decomposing $\cov\{
X^{(1)}(t), X^{(1)}(s)\}$
into eigenfunctions and eigenvalues. This leads to
$\cov\{X^{(1)}(t), X^{(1)}(s)\}=\sum_{k=1}^{\infty} \lambda_{k,1}
\phi_{k,1}(t)\phi_{k,1}(s)$ and the Karhunen--Lo\`eve representation
$ X_i^{(1)}(t) = \mu^{(1)}(t) +
\sum_{k=1}^{\infty}\xi_{ik,1}\phi_{k,1}(t)$, with orthonormal
eigenfunctions $\phi_{k,1}$ [\citet{mull091}].

%s2.2 ###
\subsection{Empirical stochastic differential equation}

In the following we consider differentiable Gaussian processes, for
which the differential equation introduced below automatically
applies. In the absence of the Gaussian assumption, one may invoke
an alternative least squares-type interpretation. Gaussianity of the
processes implies the joint normality of centered processes $\{X(t)-\mu(t),
X^{(1)}(t)-\mu^{(1)}(t)\}$ at all points $t \in\mt$, so that
%
%e5 ###
%
\begin{eqnarray} \label{jt}
&&
\pmatrix{X^{(1)}(t)-\mu^{(1)}(t)
\cr X(t)-\mu(t)}\nonumber\\
&&\qquad=\pmatrix{\displaystyle\sum_{k=1}^\infty\xi_{k}\phi^{(1)}_k(t)\cr
\displaystyle\sum_{k=1}^\infty\xi_{k}\phi_k(t)}
\\
&&\qquad\sim N_2
\left(
\pmatrix{0\cr0},
\pmatrix{
\displaystyle\sum_{k=1}^\infty\lambda_k \phi_k^{(1)}(t)^2 &
\displaystyle\sum_{k=1}^\infty
\lambda_k\phi_k^{(1)}(t) \phi_k(t)\cr
\displaystyle\sum_{k=1}^\infty\lambda_k\phi_k^{(1)}(t) \phi_k(t) &
\displaystyle\sum_{k=1}^\infty\lambda_k \phi_k(t)^2}\right).\nonumber
\end{eqnarray}

This joint normality immediately implies a ``population''
differential equation of the form $E\{X^{(1)}(t)-\mu^{(1)}(t)|
X(t)\}=\beta(t)\{X(t)-\mu(t)\}$, as\vadjust{\goodbreak} has been observed in \citet
{mull091}; for
additional details see Appendix~\ref{secA1}. However, it is considerably
more interesting to find a dynamic equation which applies to the
individual trajectories of processes $X$. This goal necessitates
inclusion of a stochastic term which leads to an empirical
stochastic differential equation that governs the dynamics of
individual trajectories $X_i$.

\begin{theorem} \label{thm1}
For a differentiable Gaussian process, it
holds that
%
%e6 ###
%
\begin{equation} \label{de}
X^{(1)}(t)-\mu^{(1)}(t)=\beta(t)\{X(t)-\mu(t)\} + Z(t),\qquad t \in\mt,
\end{equation}
where
%
%e7 ###
%
\begin{eqnarray} \label{beta}
\beta(t) &=& \frac{\cov\{X^{(1)}(t),X(t)\}}{\var\{
X(t)\}}
=\frac{\sum_{k=1}^\infty\lambda_k\phi_k^{(1)}(t)
\phi_k(t)}{\sum_{k=1}^\infty
\lambda_k\phi_k(t)^2}\nonumber\\[-8pt]\\[-8pt]
&=&\frac{1}{2}\frac{d}{dt}\log[\var\{X(t)\}],\qquad
t\in\mt,\nonumber
\end{eqnarray}
and $Z$ is a Gaussian process such that
$Z(t), X(t)$ are independent at each $t\in\mt$ and where $Z$ is
characterized by $E\{Z(t)\}=0$ and $\cov\{Z(t),Z(s)\}=G_z(t, s)$,
with
%
%e8 ###
%
\begin{eqnarray}\label{Gz}
G_z(t, s)&=&\sum_{k=1}^\infty
\lambda_k\phi_k^{(1)}(t) \phi^{(1)}_k(s) - \beta(t)
\sum_{k=1}^\infty\lambda_k\phi_k(t) \phi^{(1)}_k(s)
\nonumber\\[-8pt]\\[-8pt]
&&{} - \beta(s)\sum_{k=1}^\infty\lambda_k\phi^{(1)}_k(t)
\phi_k(s) + \beta(t)\beta(s)\sum_{k=1}^\infty\lambda_k\phi_k(t)
\phi_k(s).\nonumber
\end{eqnarray}
\end{theorem}

Equation (\ref{de}) provides a first order linear differential
equation which includes a time-varying linear coefficient function
$\beta(t)$ and a random drift process $Z(t)$. The process $Z$
``drives'' the equation at each time $t$. It is square integrable
and possesses a smooth covariance function and smooth
trajectories. It also provides an alternative characterization of
the individual trajectories of the process. The size of its
variance function $\var(Z(t))$ determines the importance of the
role of the stochastic drift component.

We note that the assumption of differentiability of the process
$X$ in Theorem~\ref{thm1} can be relaxed. It is sufficient to require weak
differentiability, assuming that $X \in W^{1,2}$, where
$H^1=W^{1,2}$ denotes the Sobolev space of square integrable
functions with square integrable weak derivative [\citet{ziem89}].
Along these lines, equation (\ref{de}) may be interpreted as a
stochastic Sobolev embedding. Observe also that the drift term $Z$
can be represented as an integrated diffusion process. Upon
combining (\ref{kl3}) and (\ref{de}), and observing that
functional principal components can be represented as
$\xi_k=\sqrt{\lambda_k/\gamma_k} \int_{\T} h_{k}(u) \,dW(u)$,
where $h_{k}$ is the $k$th eigenfunction of the Wiener process
$W$ on domain $\T=[0,T]$\vadjust{\goodbreak} and $\gamma_k$ the associated eigenvalue,
such a representation is given by
\[
Z(t)=\sum_{k=1}^{\infty} \sqrt{\frac{\lambda_k
}{2T^3}}(2k-1)\pi\int_0^T \sin\biggl\{\frac{(2k-1)\pi}{2T} u\biggr\}
\bigl\{\phi_k^{(1)}(t)-\beta(t)\phi(t)\bigr\} \,dW(u).
\]

Another observation is that the joint normality in (\ref{jt}) can be
extended to joint normality for any finite number of derivatives,
assuming these are well defined. Therefore, higher order stochastic
differential equations can be derived analogously to (\ref{de}).
However, these higher-order analogues are likely to be much less
relevant practically, as higher-order derivatives of mean and
eigenfunctions cannot be well estimated for the case of sparse noisy
data or even denser noisy data.

Finally, it is easy to see that the differential equation
(\ref{de}) is equivalent to the following stochastic integral
equation:
%
%e9 ###
%
\begin{eqnarray} \label{int}
X(t)&=&X(s)+\{\mu(t)-\mu(s)\} \nonumber\\
&&{} +
\int_s^t\beta(u)\{X(u)-\mu(u)\} \,du+\int_s^t Z(u) \,du\\
\eqntext{\mbox{for any } s,t \in\mt, 0 \le s <t,}
\end{eqnarray}
in the sense that $X$ is the solution of both equations. For
a domain with left endpoint at time $0$, setting $s=0$ in
(\ref{int}) then defines a classical initial value problem. Given
a trajectory of the drift process $Z$ and a varying coefficient
function $\beta$, one may obtain a solution for $X$ numerically by
Euler or Runga--Kutta integration or directly by applying the known
solution formula for the initial value problem of an inhomogeneous
linear differential equation.

%s2.3 ###
\subsection{Interpretations and decomposition of variance}

We note that equations (\ref{de}) and (\ref{int}) are of particular
interest on domains $\mt$ or subdomains defined by those times $t$
for which the variance function $\var\{Z(t)\}$ is ``small.'' From
(\ref{beta}) and (\ref{Gz}) one finds
%
%e10 ###
%
\begin{eqnarray}\label{v}
V(t)&=&\var\{Z(t)\}\nonumber\\
&=&\bigl({\var\bigl\{X^{(1)}(t)\bigr\}\var\{X(t)\} -
\bigl[\cov\bigl\{X^{(1)}(t),
X(t)\bigr\}\bigr]^2}\bigr)/\var\{X(t)\}\nonumber\\[-8pt]\\[-8pt]
&=& \Biggl(\sum_{k=1}^\infty\lambda_k\bigl(\phi_k^{(1)}(t)\bigr)^2
\sum_{k=1}^\infty\lambda_k \phi_k^2(t) -
\Biggl\{\sum_{k=1}^\infty\lambda_k\phi_k^{(1)}(t)
\phi_k(t) \Biggr\}^2 \Biggr)\nonumber\\
&&{}\Bigg/\sum_{k=1}^\infty
\lambda_k\phi_k^2(t).\nonumber
\end{eqnarray}
On subdomains with small variance
function, the solutions of (\ref{de}) will not deviate too much
from the solutions of the approximating equation
%
%e11 ###
%
\begin{equation}\label{approx}
X^{(1)}(t)-\mu^{(1)}(t)=\beta(t)\{X(t)-\mu(t)\},\qquad t \in\mt.\vadjust{\goodbreak}
\end{equation}
In this
situation, the future changes in value of individual trajectories
are highly predictable and the interpretations of the dynamic
behavior of processes $X$ obtained from the shape of the varying
coefficient function $\beta(t)$ apply at the individual level.

If $\beta(t)<0$, the dynamic behavior can be characterized as
``dynamic regression to the mean''; a trajectory which is away from
(above or below) the mean function $\mu$ at time $t$ is bound to
move closer toward the mean function $\mu$ as time progresses
beyond $t$. Similarly, if $\beta(t)>0$, trajectories will exhibit
``explosive'' behavior, since deviations from the mean (above or
below) at time $t$ will be reinforced further as time progresses, so
that trajectories are bound to move further and further away from
the population mean trajectory. Intermediate cases arise when the
function $\beta$ changes sign, in which case the behavior will
switch between explosive and regression to the mean, depending on
the time subdomain. Another situation occurs on subdomains where
both $\beta$ and $\var(Z(t))$ are very small, in which case the
deviation of the derivative of an individual trajectory from the
population mean derivative will also be small which means that
trajectory derivatives will closely track the population mean
derivative on such subdomains.

The independence of $Z(t)$ and $X(t)$ means that the right-hand side
of (\ref{de}) provides an orthogonal decomposition of $X^{(1)}(t)$ into
the two components $\beta(t)X(t)$ and $Z(t)$ such that
\[
\var\bigl\{X^{(1)}(t)\bigr\}=\beta(t)^2 \var\{X(t)\}+\var\{Z(t)\}.
\]
It is therefore of interest to determine the fraction of the
variance of $X^{(1)}(t)$ that is explained by the differential equation
itself, that is, the ``coefficient of determination''
%
%e12 ###
%
\begin{equation}\label{R1}
R^2(t)=\frac{\var\{\beta(t)X(t)\}}{\var\{X^{(1)}(t)\}}=1-\frac
{\var\{
Z(t)\}}{\var\{X^{(1)}(t)\}},
\end{equation}
which is seen to be equivalent to the squared correlation between
$X(t), X^{(1)}(t)$,
%
%e13 ###
%
\begin{equation}\label{R2}\qquad
R^2(t)=\frac{[\cov\{X(t),X^{(1)}(t)\}]^2}{\var\{X(t)\}\var\{
X^{(1)}(t)\}}=
\frac{\{\sum_{k=1}^\infty\lambda_k\phi_k^{(1)}(t)
\phi_k(t)\}^2}{\sum_{k=1}^\infty
\lambda_k\phi_k(t)^2 \sum_{k=1}^\infty
\lambda_k\phi_k^{(1)}(t)^2}.
\end{equation}

We are then particularly interested in subdomains of $\mt$ where
$R^2(t)$ is large, say, exceeds a prespecified threshold of 0.8 or
0.9. On such subdomains the drift process $Z$ is relatively small
compared to $X^{(1)}(t)$ so that the approximating deterministic first
order linear differential equation (\ref{approx}) can substitute for
the stochastic dynamic equation (\ref{de}). In this case,
short-term prediction of $X(t+\Delta)$ may be possible for small
$\Delta$, by
directly perusing the approximating differential
equation~(\ref{approx}).

\begin{figure}

\includegraphics{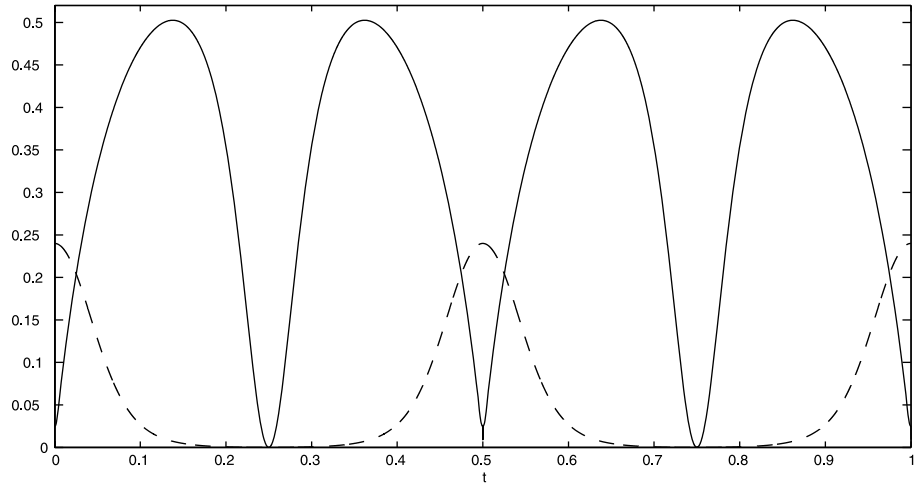}

\caption{``Coefficient of determination'' functions
$R^2(t)$ (\protect\ref{R1}), (\protect\ref{R2}), quantifying the
fraction of variance explained by the deterministic part of the dynamic
equation (\protect\ref{de}), illustrated
for the trigonometric basis $\{\sqrt{2}\cos(2k\pi t), k
=1,2,\ldots\}$ on $[0,1]$ and eigenvalue sequences $\lambda_k=k^{-4}$
(solid) and
$\lambda_k=2^{-k}$ (dashed).}\label{rfig}
\end{figure}

It is instructive to visualize an example of the function $R^2(t)$
for the case of fully specified eigenfunctions and eigenvalues.
Assuming\vspace*{1pt} that the eigenfunctions correspond\vadjust{\goodbreak} to the trigonometric
orthonormal system $\{\sqrt{2}\cos(2k\pi t), k =1,2,\ldots\}$ on
$[0,1]$, we find from (\ref{R2})
\begin{eqnarray*}
R^2(t)&=& \Bigl[\sum\lambda_k k
\cos(2k\pi t)\sin(2k\pi t) \Bigr]^2\\
&&{}\big/ \Bigl[\sum\lambda_k (\cos
(2k\pi
t))^2 \sum\lambda_k k (\sin(2k\pi t))^2 \Bigr],\qquad t \in
[0,1].
\end{eqnarray*}
Choosing $\lambda_k=k^{-4}, \lambda_k=2^{-k}$ and
numerically approximating these sums, one obtains the functions
$R^2(t)$ as depicted in Figure~\ref{rfig}. This illustration shows
%
%f1 ###
%
that the behavior of this function often will fluctuate between
small and large values and also depends critically on both the
eigenvalues and the shape of the eigenfunctions.

%s3 ###
\section{Asymptotic properties}\label{sec3}

We obtain asymptotic consistency results for estimators of the
varying coefficient functions $\beta$, for the variance function
$\var\{Z(t)\}$ of the drift process and for the variance explained at
time $t$ by the deterministic part (\ref{approx}) of the
stochastic equation (\ref{de}), quantified by $R^2(t)$.
Corresponding estimators result from plugging in estimators for the
eigenvalues $\lambda_k$, eigenfunctions $\phi_k$ and eigenfunction
derivatives $\phi_k^{(1)}$ into the representations (\ref{beta}) for
the function $\beta(t)$, (\ref{v}) for the variance function of $Z$
and (\ref{R2}) for $R^2(t)$. Here one needs to truncate the
expansions at a finite number $K=K(n)$ of included eigen-components.

Details about the estimation procedures, which are based on local
linear smoothing of one- and two-dimensional functions, are
deferred to Appendix~\ref{secA2}. Our asymptotic consistency results focus
on $L^2$ convergence rates. They peruse auxiliary results on the
convergence of estimates of eigenvalues, eigenfunctions and
eigenfunction derivatives, complementing and improving upon related
results of \citet{mull091}, which were derived for convergence
in the
sup norm. Improved rates of convergence in the $L^2$ distance are
the consequence of a special decomposition that we employ in the
proofs to overcome the difficulty caused by the dependence of the
repeated measurements.

Required regularity conditions include assumptions for the
distribution of the design points, behavior of eigenfunctions
$\phi_k$ and eigenvalues $\lambda_k$ as their order $k$ increases
and the large sample behavior of the bandwidths $h_{\mu,0},
h_{\mu,1}$ for the estimation of the mean function $\mu$ and its
first derivative $\mu^{(1)}(t)$, and $h_{G,0}, h_{G,1}$ for the
estimation of the covariance surface and its partial derivative. We
note that extremely sparse designs are covered, with only two
measurements per trajectory; besides being bounded, the number of
measurements $N_i$ for the $i$th trajectory is required to satisfy
$P(N_i \ge2)>0$.

Specifically, for the observations $(T_{ij}, Y_{ij})$, $i =
1,\ldots, n$, $j = 1,\ldots,N_i$, made for the $i$th trajectory, we
require that:
\begin{enumerate}[(A1)]%\label{a1}
\item[(A1)] $N_i$ are random variables with $N_i \stackrel{\mathrm
{i.i.d.}}{\sim} N$,
where $N$ is a bounded positive discrete random variable with
and $P\{N \ge2\} > 0$, and $(\{T_{ij}, j \in J_i \}, \{Y_{ij}, j \in
J_i\})$ are independent of $N_i$,
for $J_i \subseteq\{1, \ldots, N_i\}$.
\end{enumerate}
Writing $\mathbf{T}_i = (T_{i1}, \ldots, T_{iN_i})^T$ and $\mathbf
{Y}_i = (Y_{i1}, \ldots, Y_{iN_i})^T$, the triples $\{\mathbf{T}_i,
\mathbf{Y}_i, N_i\}$ are assumed to be i.i.d. For the bandwidths used
in the
smoothing steps for $\mu(t)$ and $\mu^{(1)}(t)$ in (\ref{smooth1}),
$G(t,s)$ and $G^{(1,0)}(t,s)$ in (\ref{smooth2}), we require that,
as $n\rightarrow\infty$,
\begin{enumerate}[(A2)]
\item[(A2)] $\max(h_{\mu,0}, h_{\mu,1}, h_{G,0}, h_{G,1})
\rightarrow0$, $nh_{\mu,0} \rightarrow\infty$, $nh_{\mu,1}^3
\rightarrow\infty$,
$nh_{G,0}^2 \rightarrow\infty$, $nh_{G,1}^4 \rightarrow\infty$.
\end{enumerate}
To characterize the behavior of estimated eigenfunction
derivatives $\hat{\phi}^{(1)}(t)$, define
%
%e14 ###
%
\begin{equation}\label{spacing}
\delta_1=\lambda_1-\lambda_2, \delta_k=\min_{j\leq k}
(\lambda_{j-1}-\lambda_j, \lambda_j-\lambda_{j+1}),\qquad k\geq2.
\end{equation}
For the kernels used in the local linear smoothing steps and
underlying density and moment functions, we require assumptions
(B1) and (B2) in the \hyperref[app]{Appendix}. Denote the $L^2$ norm by
$\|f\|=\{\int_\T f^2(t) \,dt\}^{1/2}$, the Hilbert--Schmidt norm by
$\|\Phi\|_s=\{\int_\T\int_\T\{\Phi^2(t, s) \,dt \,ds\}^{1/2}$ and
also define $\|\Phi\|_u^2=\{\int_\T\Phi^2(t, t)
\,dt\}^{1/2}$.\vspace*{1pt}

The following result provides asymptotic rates of convergence in the
$L^2$ norm for the auxiliary estimates of mean functions and their
derivatives as well as covariance functions and their partial
derivatives, which are briefly discussed in Appendix~\ref{secA2}. A
consequence is a convergence result for the eigenfunction derivative
estimates $\hat{\phi}_k^{(1)}$, with constants and rates that hold
uniformly in the order $k \ge1$.
\begin{theorem}\label{thm2}
Under \textup{(A1)} and \textup{(A2)} and \textup{(B1)--(B3)},
for $\nu\in\{0, 1\}$,
%
%e15 ###
%
\begin{eqnarray}
\label{thm2-eq1}
\bigl\|\hat{\mu}^{(\nu)}-\mu^{(\nu)}\bigr\| &=& O_p \biggl(\frac{1}{\sqrt
{nh_{\mu,\nu}^{2\nu+1}}}+h_{\mu,\nu}^2 \biggr),\nonumber\\[-8pt]\\[-8pt]
\bigl\|\hat{G}^{(\nu,0)}-G^{(\nu,0)}\bigr\|_s &=& O_p \biggl(\frac{1}{\sqrt
{n}h_{G,\nu}^{\nu+1}}+h_{G,\nu}^2 \biggr).\nonumber
\end{eqnarray}
For $\phi^{(1)}_k(t)$ corresponding to $\lambda_k$ of
multiplicity $1$,
%
%e16 ###
%
\begin{eqnarray}\label{thm2-eq2}
&&
\bigl\|\hat{\phi}^{(1)}_k(t)-\phi^{(1)}_k(t)\bigr\|\nonumber\\[-8pt]\\[-8pt]
&&\qquad=O_p \biggl(\frac
{1}{\lambda_k}
\biggl\{\frac{1}{\sqrt{n}h_{G,
1}^{2}}+h_{G,1}^2+\frac{1}{\delta_k} \biggl(\frac{1}{\sqrt{n}h_{G,
0}}+h_{G,0}^2 \biggr) \biggr\} \biggr),\nonumber
\end{eqnarray}
where the $O_p(\cdot)$ term in
(\ref{thm2-eq2}) is uniform in $k \ge1$.
\end{theorem}

An additional requirement is that variances of processes $X$ and
$X^{(1)}$ are bounded above and below, since these appear in the
denominators of various representations, for example, in (\ref{v}) and
(\ref{R2}),
\begin{enumerate}[(A3)]
\item[(A3)] $\inf_{t \in\T} G^{(\nu, \nu)}(s, s)\geq c>0$ and $\|
G^{(\nu, \nu)}\|_u<\infty$ for $\nu=0,
1$,
\end{enumerate}
implying that $\|G^{(\nu, \nu)}\|_s<\infty$ by the
Cauchy--Schwarz inequality. Define remainder terms
%
%e17 ###
%
\begin{equation} \label{rem}
R_{K,\nu}(t)=\sum_{k=K+1}^\infty\lambda_k\bigl\{\phi^{(\nu)}(t)\bigr\}^2,\qquad
R^\ast_{K, \nu}(s, t)=\sum_{k=K+1}^\infty
\lambda_k\phi^{(\nu)}(s)\phi^{(\nu)}(t);\hspace*{-28pt}
\end{equation}
by the Cauchy--Schwarz inequality, $\|R^\ast_{K, \nu}\|_s\leq
\|R_{K,\nu}\|_u$.

In order to obtain consistent estimates of various quantities, a
necessary requirement is that the first $K$ eigen-terms approximate
the infinite-dimensional process sufficiently well. The increase in
the sequence $K=K(n)$ as $n \rightarrow\infty$ therefore needs to
be tied to the spacing and decay of eigenvalues,
\begin{enumerate}[(A4)]
\item[(A4)]
\begin{eqnarray*}\\[-35pt]
K&=&o\bigl(\min\bigl\{\sqrt{n}h_{G,1}^{2}, h_{G,1}^{-2}\bigr\}\bigr),\\
\sum_{k=1}^K \delta_k^{-1}&=&o\bigl(\min\bigl(\sqrt{n}h_{G,0},
h_{G,0}^{-2}\bigr\}\bigr),\\
\max_{\nu=0, 1}\|R_{K, \nu}\|&\rightarrow&0\qquad\mbox{as $n\rightarrow
\infty$.}
\end{eqnarray*}
\end{enumerate}
If the eigenvalues decrease rapidly and merely a few leading
terms are needed, condition (A4) is easily satisfied. We use
``$\stackrel{p}{\asymp}$'' to connect two terms which are
asymptotically of the same\vadjust{\goodbreak} order in probability, that is, the terms are
$O_p$ of each other. Define the sequence
%
%e18 ###
%
\begin{equation}\label{rate}\quad
\alpha_n=K\bigl\{\bigl(\sqrt{n}h_{G,1}^{2}\bigr)^{-1}+h_{G,1}^2\bigr\}+\Biggl(\sum_{k=1}^K
\delta_k^{-1}\Biggr)\bigl\{\bigl(\sqrt{n}h_{G,0}\bigr)^{-1}+h_{G,0}^2\bigr\}.
\end{equation}

Note that $\mbox{cov}\{X^{(1)}(s), X^{(1)}(t)\}=G^{(1,1)}(s,
t)=\sum_{k=1}^\infty\lambda_k \phi_k^{(1)}(s)\phi_k^{(1)}(t)$ with
corresponding plug-in estimate $\hat{G}_K^{(1,1)}(s, t)=\sum_{k=1}^K
\hat{\lambda}_k \hat{\phi}_k^{(1)}(s)\hat{\phi}_k^{(1)}(t)$, where
$K=K(n)$ is the included number of eigenfunctions. The plug-in
estimate for $\beta(t)$ is based on (\ref{beta}) and given by
$\hat{\beta}_K(t)=\sum_{k=1}^K \hat{\lambda}_k\hat{\phi}_k^{(1)}(t)
\hat{\phi}_k(t)/\sum_{k=1}^K \hat{\lambda}_k\hat{\phi}_k(t)^2$ and
analogously the plug-in estimate $\hat{G}_{z,K}$ of $G_z$ is based
on representation~(\ref{Gz}), using the estimate $\hat{\beta}_K$.
In a completely analogous fashion one obtains the estimates
$\hat{R}_K^2(t)$ of $R^2(t)$ from (\ref{R2}) and $\hat{V}_K(t)$ of
the variance function $V(t)=\var(Z(t))$ of the drift process from
(\ref{v}). The $L^2$ convergence rates of these estimators of
various components of the dynamic model (\ref{de}) are given in the
following result.
\begin{theorem}\label{thm3}
Under \textup{(A1)--(A4)} and \textup{(B1)--(B3)},
%
%e20 ###
%e19 ###
%
\begin{eqnarray}
\label{thm3-eq1}
\bigl\|\hat{G}_K^{(1,1)}-G^{(1,1)}\bigr\|_s&=&O_p(\alpha_n+\|R_{K,1}^\ast\|
_s),\nonumber\\[-8pt]\\[-8pt]
\bigl\|\hat{G}_K^{(1,1)}-G^{(1,1)}\bigr\|_u&=&O_p(\alpha_n+\|R_{K,1}\|),
\nonumber\\
\|\hat{\beta}_K-\beta\| &\stackrel{p}{\asymp}& \|\hat{G}_{z,K}-
G_z\|_s \stackrel{p}{\asymp} \|\hat{G}_{z,K}-G_z\|_u\nonumber\\
\label{thm3-eq2}
&\stackrel{p}{\asymp}& \|\hat{R}^2_K-R^2\|
\stackrel{p}{\asymp} \|\hat{V}_K-V\| \\
&=&O_p(\alpha_n+\|R_{K,0}\|+\|R_{K,1}\|).\nonumber
\end{eqnarray}
\end{theorem}

The weak convergence and $L^2$ consistency for the estimated
eigenvalues $\{\rho_k\}$ and eigenfunctions $\{\psi_k\}$ of the
drift process $Z$ is an immediate consequence of this result. To see
this,\vspace*{1pt} one may use $\sup_{k\geq1}
|\hat{\rho}_k-\rho_k|=O_p(\|\hat{G}_{z}-G_z\|_s)$ and
$\|\hat{\psi}_k-\psi_k\|=\delta_k^{\ast-1}
O_p(\|\hat{G}_{z}-G_z\|_s)$ where $\hat{G}_{z}$ is any estimate of
$G_z$ [\citet{bosq00}]. Here the $O_p(\cdot)$ terms are uniform
in $k$
and $\delta_1^\ast=\rho_1-\rho_2, \delta_k^\ast=\min
_{j\leq k}
(\rho_{j-1}-\rho_j, \rho_j-\rho_{j+1})$ for $k\geq2$.

%s4 ###
\section{Application to online auction data}\label{sec4}

%s4.1 ###
\subsection{Data and population level analysis}\label{sec41}

To illustrate our methods, we analyze the dynamic system
corresponding to online auction data, specifically using eBay
bidding data for 156 online auctions of Palm Personal Digital
Assistants in 2003 (courtesy of Wolfgang Jank). The data are
``live bids'' that are entered by bidders at irregular times and
correspond to the actual price a winning bidder would pay for the
item. This price is usually lower than the ``willingness-to-pay''
price, which is the value a bidder enters. Further details
regarding the proxy bidding mechanism for the 7-day second-price
auction design that applies to these data can be found in
Jank and Shmueli (\citeyear{jank051}, \citeyear{jank062}),
Liu and M\"{u}ller (\citeyear{mull081}, \citeyear{mull091}).

The time unit of these 7-day auctions is hours and the domain is the
interval $[0,168]$. Adopting the customary approach, the bid prices
are log-transformed prior to the analysis. The values of the live
bids $Y_{ij}$ are sampled at bid arrival times $T_{ij}$, where
$i=1,\ldots, 156$ refers to the auction index and $j=1,\ldots,N_i$
to the total number of bids submitted during the $i$th auction; the
number of bids per auction is found to be between 6 and 49 for these
data. We adopt the point of view that the observed bid prices result
from an underlying price process which is smooth, where the bids
themselves are subject to small random aberrations around underlying
continuous trajectories. Since there is substantial variability of
little interest in both bids and price curves during the first three
days of an auction, when bid prices start to increase rapidly from a
very low starting point to more realistic levels, we restrict our
analysis to the interval [$96, 168$] (in hours), thus omitting the
first three days of bidding. This allows us to focus on the more
interesting dynamics in the price curves taking place during the
last four days of these auctions.

Our aim is to explore the price dynamics through the empirical
stochastic differential equation (\ref{de}). Our study emphasizes
description of the dynamics over prediction of future auction
prices and consists of two parts: a description of the dynamics of
the price process at the ``population level'' which focuses on
patterns and trends in the population average and is reflected by
dynamic equations for conditional expectations. The second and
major results concern the quantification of the dynamics of
auctions at the individual or ``auction-specific level'' where
one studies the dynamic behavior for each auction separately, but
uses the information gained across the entire sample of auctions.
Only the latter analysis involves the stochastic drift term $Z$ in
the stochastic differential equation (\ref{de}). We begin by
reviewing the population level analysis, which is characterized by
the deterministic part of (\ref{de}), corresponding to the
equation $E(X^{(1)}(t)-\mu^{(1)}(t)|X(t)-\mu(t))=\beta(t)\{X(t)-\mu
(t)\}$. This equation describes a
relationship that holds for conditional means but not necessarily
for individual trajectories.

For the population level analysis, we require estimates of the mean
price curve $\mu$ and its first derivative $\mu^{(1)}$, and these
are obtained by applying linear smoothers to (\ref{smooth1}) to the
pooled scatterplots that are displayed in Figure~\ref{auc-mu} (for
more details, see Appendix~\ref{secA2}). One finds that both log prices and
log price derivatives are increasing throughout, so that at the
log-scale the price increases are accelerating in the mean as the
auctions proceed.

%f2 ###
%
\begin{figure}

\includegraphics{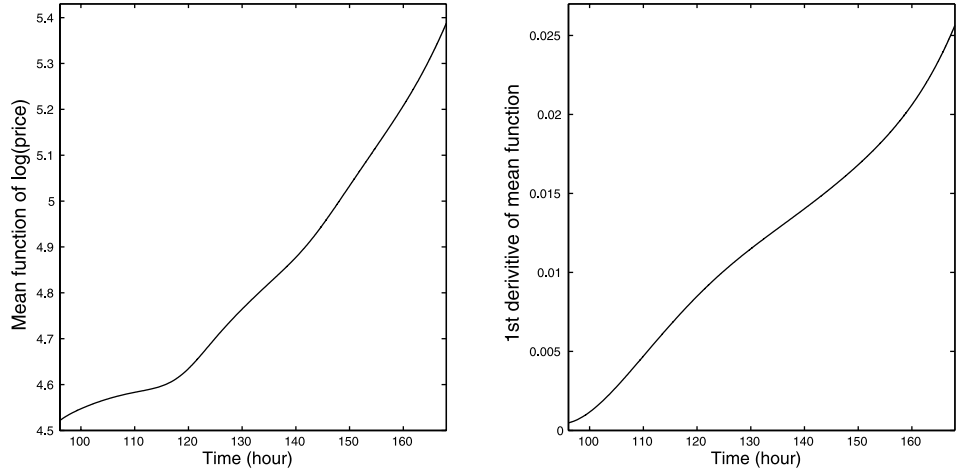}

\caption{Smooth estimate of the mean function of
log(Price) in the left panel and of its first derivative in the
right panel.}\label{auc-mu}
\end{figure}

%f3 ###
%
\begin{figure}[b]

\includegraphics{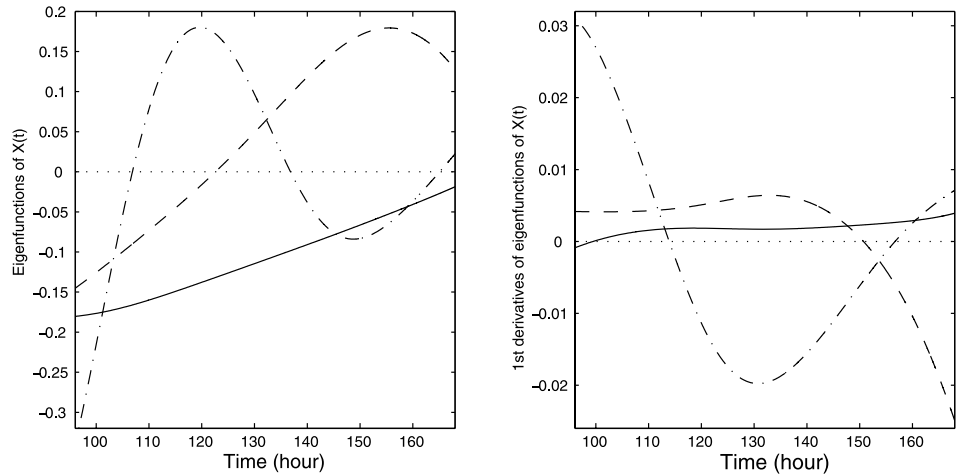}

\caption{Smooth estimates of the first (solid),
second (dashed) and third (dotted) eigenfunctions of process $X$
(left panel) and of their derivatives (right panel),
$\hat{\phi}_1^{(1)}$ (solid), $\hat{\phi}_2^{(1)}$ (dashed) and
$\hat{\phi}_3^{(1)}$ (dash-dotted).}\label{auc-xeig}
\end{figure}

%using the bandwidths $h_{\mu,0}=7, h_{\mu, 1}=15, h_{G,0}= h_{G,
%1}=15$, respectively.
A second ingredient for our analysis are estimates for the
eigenfunctions and eigenvalues (details in Appendix~\ref{secA2}). Since
the first three eigenfunctions were found to explain 84.3\%,
14.6\% and 1.1\% of the total variance, three components were
selected. The eigenfunction estimates are shown in the left panel
of Figure~\ref{auc-xeig}, along with the estimates of the
corresponding eigenfunction derivatives in the right panel. For
the interpretation of the eigenfunctions it is helpful to note
that the sign of the eigenfunctions is arbitrary. We also note
that variation in the direction of the first eigenfunction
$\phi_1$ corresponds to the major part of the variance. The
variances $\lambda_1 \phi_1^2(t)$ that are attributable to this
eigenfunction are seen to steadily decrease as $t$ is increasing,
so that this eigenfunction represents a strong trend of higher
earlier and smaller later variance in the log price trajectories.

The contrast between large variance of the trajectories at earlier
times and smaller variances later reflects the fact that auction
price trajectories are less determined early on when both
relatively high as well as low prices are observed, while at later
stages prices differ less as the end of the auction is approached
and prices are constrained into a narrower range. Correspondingly,
the first eigenfunction derivative is steadily increasing
(decreasing if the sign is switched), with notably larger
increases (decreases) both at the beginning and at the end and a
relatively flat positive plateau in the middle part.

The second eigenfunction corresponds to a contrast between
trajectory levels during the earlier and the later part of the
domain, as is indicated by its steady increase and the sign
change, followed by a slight decrease at the very end. This
component thus reflects a negative correlation between early and
late log price levels. The corresponding derivative is positive
and flat, with a decline and negativity toward the right
endpoint. The third eigenfunction, explaining only a small
fraction of the overall variance, reflects a more complex
contrast between early and late phases on one hand and a middle
period on the other, with equally more complex behavior reflected
in the first derivative.

The eigenfunctions and their derivatives in conjunction with the
eigenvalues determine the varying coefficient function $\beta$,
according to (\ref{beta}). The estimate of this function is
obtained by plugging in the estimates for these quantities and is
visualized in the left panel of Figure~\ref{auc-beta-zeig},
demonstrating small negative values for the function $\beta$
throughout most of the domain, with a sharp dip of the function
into the negative realm near the right end of the auctions.

For subdomains of functional data, where the varying coefficient
or ``dynamic transfer'' function $\beta$ is negative, as is the
case for the auction data throughout the entire time domain, one
may interpret the population equation $E(X^{(1)}(t)-\mu
^{(1)}(t)|X(t)-\mu(t))=\beta(t)\{X(t)-\mu(t)\}$
as indicating ``dynamic regression to the mean.'' By this we mean
the following: when a trajectory value at a current time $t$ falls
above (resp., below) the population mean trajectory value at $t$,
then the conditional mean derivative of the trajectory at $t$
falls below (resp., above) the mean. The overall effect of this
negative association is that the direction of the derivative is
such that trajectories tend to move toward the overall population
mean trajectory as time progresses.

%f4 ###
%
\begin{figure}[b]

\includegraphics{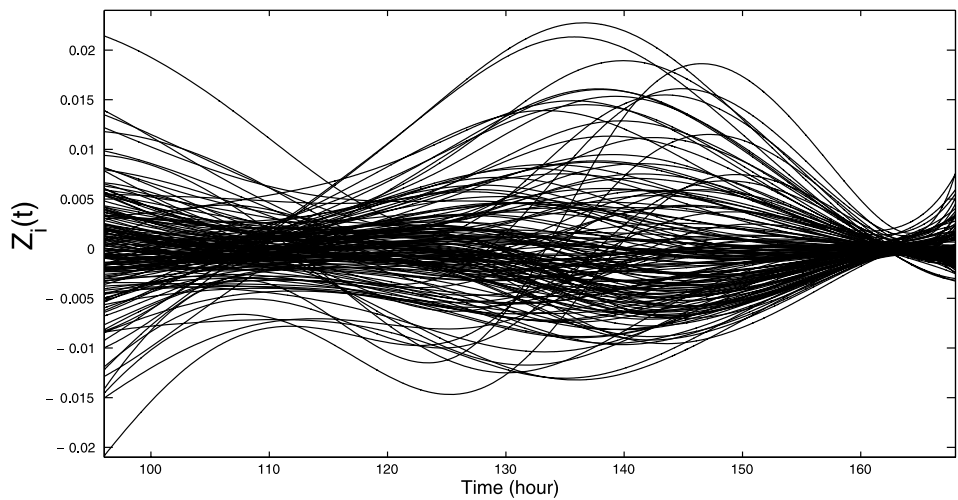}

\caption{Smooth estimates for the trajectories of the
drift process $Z$.}\label{auc-z}
\end{figure}

Thus, our findings for the auction data indicate that ``dynamic
regression to the mean'' takes place to a small extent throughout
the auction period and to a larger extent near the right tail, at
the time when the final auction price is determined [see
also \citet{mull091}]. One interpretation is that at the population
level, prices are self-stabilizing, which tends to prevent price
trajectories running away toward levels way above or below the
mean trajectory. This self-stabilization feature gets stronger
toward the end of the auction, where the actual ``value'' of the
item that is being auctioned serves as a strong homogenizing
influence. This means that in a situation where the current price
level appears particularly attractive, the expectation is that the
current price derivative is much higher than for an auction with
an unattractive (from the perspective
of a buyer) current price, for which then the corresponding
current price derivative is likely lower. The net effect is a
trend for log price trajectories to regress to the mean trajectory
as time progresses.

%s4.2 ###
\subsection{Auction-specific dynamics}

We illustrate here the proposed stochastic differential equation
(\ref{de}). First estimating the function $\beta$, we obtain the
trajectories $Z_i$ of the drift process. These trajectories are
presented in Figure~\ref{auc-z} for the entire sample of auctions.
They quantify the component of the derivative process
$X^{(1)}$ that is left unexplained by the varying coefficient
function and linear part of the dynamic model (\ref{de}). The
trajectories $Z_i$ exhibit fluctuating variances across various
subdomains. The subdomains for which these variances are small are
those where the deterministic approximation (\ref{approx}) to the
stochastic differential equation works best. It is noteworthy that
the variance is particularly small on the subdomain starting at
around 158 hours toward the endpoint of the auction at 168~hours,
since auction dynamics are of most interest during these last hours.
It is well known that toward the end of the auctions, intensive
bidding takes place, in some cases referred to as ``bid sniping,''
where bidders work each other into a frenzy to outbid each other in
order to secure the item that is auctioned.

%f5 ###
%
\begin{figure}

\includegraphics{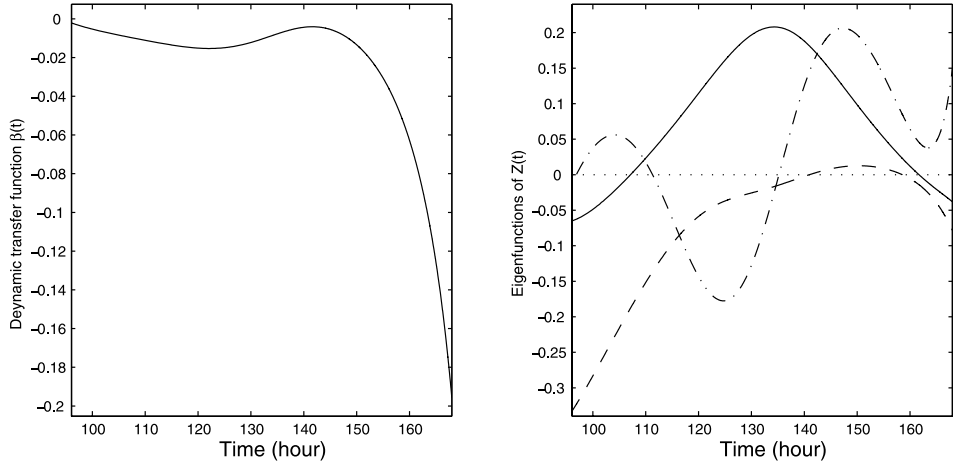}

\caption{Left: smooth estimate
of the dynamic varying coefficient function $\beta$. Right: smooth
estimates of the first (solid), second (dashed) and third
(dash-dotted) eigenfunction of $Z$ based on (\protect\ref{Gz}).}
\label{auc-beta-zeig}
\end{figure}

The right panel of Figure~\ref{auc-beta-zeig} shows the first three
eigenfunctions of $Z$, which are derived from the eigenequations
derived from estimates $\hat{G}_{z,K}$ of covariance kernels
${G}_{z,K}$ (\ref{Gz}) that are obtained as described after
(\ref{smooth2}). In accordance with the visual impression of the
trajectories of $Z$ in Figure~\ref{auc-z}, the first eigenfunction
reflects maximum variance in the middle portion of the domain and
very low variance at both ends. Interestingly, the second
eigenfunction reflects high variance at the left end of the domain
where prices are still moving upward quite rapidly, and very low
variance near the end of the auction. This confirms that overall
variation is large in the middle portion of the auctions, so that
the drift process in (\ref{de}) plays an important role in that
period.

\begin{figure}

\includegraphics{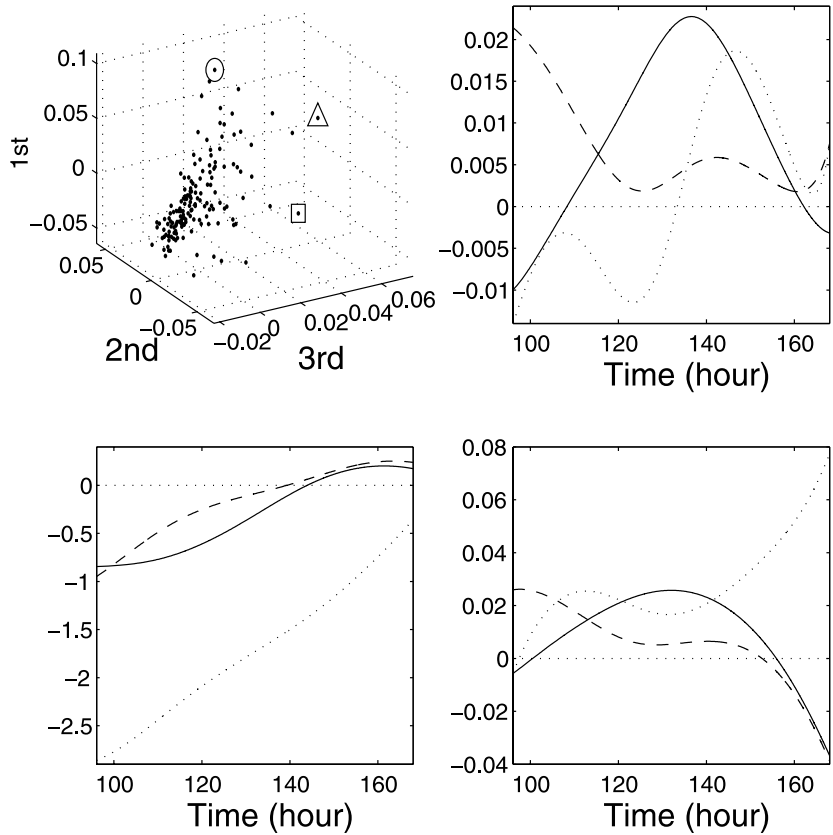}

\caption{Top left: point
cloud corresponding to the three leading FPC scores of trajectories
$Z_i$, where the point marked by a circle corresponds to the auction
with the largest (in absolute value) first score, the point marked
with a square to the auction with the largest second score and the
point marked with a ``triangle'' to the auction with the largest
third score, respectively. Top right: the trajectories $Z_i$ of the
drift process for these three auctions, where the solid curve
corresponds to the trajectory of the ``circle'' auction, the dashed
curve to the ``square'' auction and the dash-dotted curve to the
``triangle'' auction. Bottom left: corresponding centered
trajectories $X_i$. Bottom right: corresponding centered trajectory
derivatives $X_i^{(1)}$.}\label{auc-zscore-ext}
\end{figure}

Further explorations of the modes of variation of the drift process
$Z$ can be based on the functional principal component scores of
$Z_i$. Following \citet{jone92}, we identify the three auctions with
the largest absolute values of the scores. A scatterplot of second
and first principal component scores with these auctions highlighted
can be seen in the left upper panel of Figure~\ref{auc-zscore-ext}.
The corresponding individual (centered) trajectories of the drift
process $Z$ are in the\vspace*{1pt} right upper panel, and the corresponding
trajectories of centered processes $X$ and $X^{(1)}$ in the left and
right lower panels. The highlighted trajectories of $Z$ are indeed
similar to the corresponding eigenfunctions (up to sign changes), and
we find that they all exhibit the typical features of small variance
near the end of the auction for $Z$ and $X$ and of large variance
for $X^{(1)}$.

%f6 ###
%

For the two trajectories corresponding to maximal scores for first and
second eigenfunction of $Z$ we find that near the end of the auctions
their centered derivatives turn negative. This is in line with dynamic
regression to the mean, or equivalently, negative varying coefficient
function $\beta$, as described in Section~\ref{sec41}. Here the
trajectories for $X$ at a current time $t$ are above the mean
trajectory, which means the item is pricier than the average price at
$t$. As predicted by dynamic regression to the mean, log price
derivative trajectories at $t$ are indeed seen to be below the mean
derivative trajectories at $t$. The trajectory corresponding to maximal
score for the third eigenfunction also follows dynamic regression to
the mean: here the trajectory for $X$ is below the overall mean
trajectory, so that the negative varying coefficient function $\beta$
predicts that the derivative trajectory $X^{(1)}$ should be above the
mean, which indeed is the case.

%f7 ###
%
\begin{figure}

\includegraphics{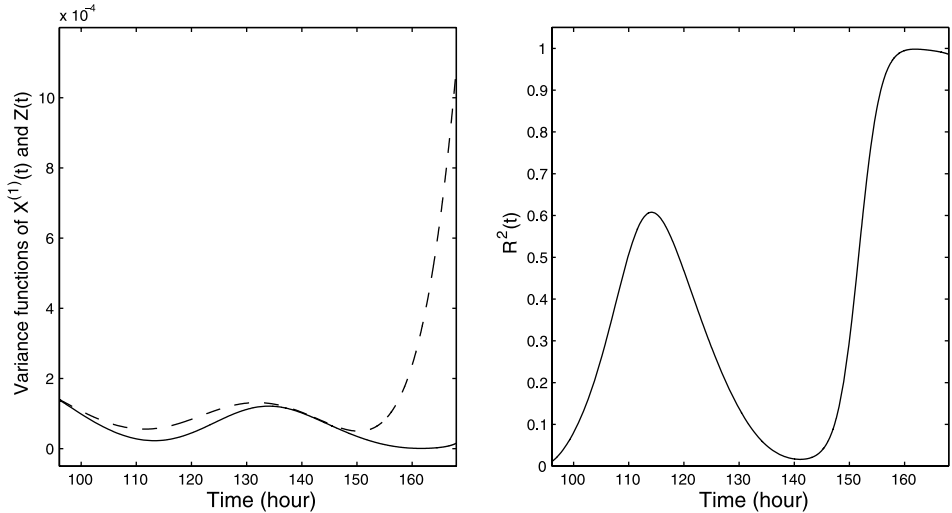}

\caption{Left: smooth estimates of the variance
functions of $X^{(1)}(t)$ (dashed) and $Z(t)$ (solid). Right: smooth
estimate of $R^2(t)$ (\protect\ref{R1}), the variance explained by the
deterministic part of the dynamic equation at time $t$.}\label{auc-var-r2}
\end{figure}

That the variance of the drift process $Z$ is small near the
endpoint of the auction is also evident from the estimated variance
function $V(t)=\var(Z(t))$ in the left panel of Figure
\ref{auc-var-r2}, overlaid with the estimated variance function
$\var(X^{(1)}(t))$ of~$X^{(1)}$. The latter is rapidly increasing
toward the end of the auction, indicating that the variance of the
derivative process is very large near the auction's end. This means
that price increases vary substantially near the end across
auctions. The large variances of derivatives coupled with the fact
that $\var(Z(t))$ is small near the end of the auction implies that
the deterministic part (\ref{approx}) of the empirical
differential equation (\ref{de}) explains a very high fraction of
the variance in the data. This corresponds to a very high, indeed
close to the upper bound 1, value of the coefficient of
determination $R^2(t)$ (\ref{R1}), (\ref{R2}) in an interval of
about 10 hours before the endpoint of an auction, as seen in the
right panel of Figure~\ref{auc-var-r2}. We therefore find that the
dynamics during the endrun of an auction can be adequately modeled by
the simple deterministic approximation (\ref{approx}) to the
stochastic dynamic equation (\ref{de}), which always applies.

\begin{figure}

\includegraphics{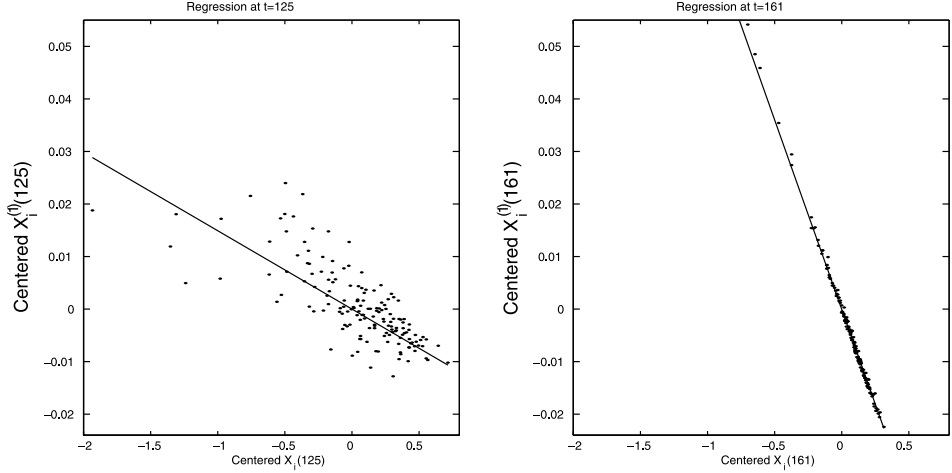}\vspace*{-3pt}

\caption{Regression of $X_i^{(1)}(t)$ on $X_i(t)$
(both centered) at $t=125$ hours (left panel) and $t=161$ hours
(right panel),\vspace*{1pt} respectively, with regression slopes
$\beta(125)=-0.015$ and coefficient of determination
$R^2(125)=0.28$, respectively, $\beta(161)=-0.072$ and
$R^2(161)=0.99$, demonstrating that the deterministic part
(\protect\ref{approx}) of the empirical differential equation
(\protect\ref{de})
explains almost the entire variance of $X^{(1)}$ at $t=161$ hours
but only a fraction of variance at $t=125$ hours.}
\label{auc-reg}\vspace*{-3pt}
\end{figure}

This finding is corroborated by visualizing the regressions of
$X^{(1)}(t)$ versus $X(t)$ at various fixed times $t$. These
regressions are linear in the Gaussian case and may be
approximated by a linear regression in the least squares sense in
the non-Gaussian case. The scatterplots of
$\hat{X}_i^{(1)}(t)-\hat{\mu}^{(1)}(t)$ versus
$\hat{X}_i(t)-\hat{\mu}(t)$ for times $t=125$ hours and $t=161$
hours (where the time domain of the auctions is between 0 and 168
hours) are displayed in Figure~\ref{auc-reg}. This reveals the
relationships to be indeed very close to linear. These are
regressions through the origin. The regression slope parameters are
not estimated from these scatterplot data which are contaminated
by noise, but rather are obtained directly from (\ref{de}), as
they correspond to $\beta(t)$. Thus one simply may use the already
available slope estimates, $\hat{\beta}(125)=-0.015$ and
$\hat{\beta}(1)=-0.072$. The associated coefficients of
determination, also directly estimated via (\ref{R2}) and the
corresponding estimation procedure,\vspace*{1pt} are found to be
$\hat{R}^2(125)=0.28$ and $\hat{R}^2(161)=0.99$.

%f8 ###
%

As the regression line fitted near the end of the auction at
$t=161$ hours explains almost all the variance, the approximating
deterministic differential equation (\ref{approx}) can be
assumed to hold at that time (and at later times as well, all the
way to the end of the auction). At $t=125$ the regression line
explains only a fraction of the variance, while a sizable portion
of variance resides in the drift process $Z$, so that the
stochastic part in the dynamic system (\ref{de}) cannot be
comfortably ignored in this time range. These relationships can be
used to predict\vadjust{\goodbreak} derivatives of trajectories and thus price changes
at time $t$ for individual auctions, given their log price
trajectory values at $t$. We note that such predictions apply to
fitted trajectories, not for the actually observed prices which
contain an additional random component that is unpredictable,
according to model (\ref{kl2}). We find that at time $t=161$,
regression to the mean is observed at the level of individual
auctions: an above (below) average log price level is closely
associated with a below (above) average log price derivative. This
implies that a seemingly very good (or bad) deal tends to be not quite
so good (or bad) when the auction ends.

%Note that bandwidths $h_{G,0}=14.5$ and $h_{G, 1}=18.2$ for $\hat{G}$
%and $\hat{G}^{(1, 0)}$.

%s5 ###
\section{Discussion}\label{sec5}

The main motivation of using the dynamic system approach based on
(\ref{de}) is that it provides a better description of the
mechanisms that drive longitudinal data but are not directly
observable. The empirical dynamic equation may also suggest
constraints on the form of parametric differential equations that
are compatible with the data. In the auction example, the dynamic
equation quantifies both the nature and extent of how expected
price increases depend on auction stage and current price level.
This approach is primarily phenomenological and does not directly
lend itself to the task of predicting future values of individual
trajectories.

That expected conditional trajectory derivatives satisfy a first-order
differential equation model (which we refer to as the
``population level'' since this statement is about conditional
expectations) simply follows from Gaussianity and in particular does
not require additional assumptions. This suffices to infer the
stochastic differential equation described in (\ref{jt}) which we term
``empirical differential equation'' as it is determined by the
data. Then the function $R^2$, quantifying the relative contribution
of the drift process $Z$ to the variance of $X^{(1)}$, determines
how closely individual trajectories follow the deterministic part of
the equation. We could equally consider stochastic differential
equations of other orders, but practical considerations favor the
modeling with first-order equations.

We find in the application example that online auctions follow a
dynamic regression to the mean regime for the entire time domain,
which becomes more acute near the end of the auction. This allows us
to construct predictions of log price trajectory derivatives from
trajectory levels at the same~$t$. These predictions get better
toward the right endpoint of the auctions. This provides a
cautionary message to bidders, since
an auction that looks particularly promising since it has a
current low log price trajectory
is likely not to stay that way and larger than average
price increases are expected down the line. Conversely, an auction
with a seemingly above average log price trajectory is likely found to have
smaller than average price increases down the line.

This suggests
that bidders take a somewhat detached stance, watching auctions
patiently as they evolve. In particular, discarding
auctions that appear overpriced is likely not a good strategy as
further price increases are going to be\vadjust{\goodbreak} smaller than the average
for such auctions. It also implies that bid snipers are ill
advised: a seemingly good deal is not likely to stay that way,
suggesting a more relaxed stance. Conversely, a seller who
anxiously follows the price development of an item, need not
despair if the price seems too low at a time before closing, as it
is likely to increase rapidly toward the end of the auction.

For prediction purposes, drift processes $Z_i$ for individual
auctions are of great interest. In time domains where their
variance is large, any log price development is possible.
Interestingly, the variance of drift processes is very small
toward the right tail of the auctions, which means that the
deterministic part of the differential equation (\ref{de}) is
relatively more important, and log price derivatives during the
final period of an auction become nearly deterministic and thus
predictable.

Other current approaches of statistical modeling of differential
equations for time course data [e.g., \citet{rams07}] share the
idea of modeling with a first order equation. In all other regards
these approaches are quite different, as they are based on the
prior notion that a differential equation of a particular and
known form pertains to the observed time courses and moreover
usually have been developed for the modeling of single time
courses. This established methodology does not take into account
the covariance structure of the underlying stochastic process. In
contrast, this covariance structure is a central object in our
approach and is estimated nonparametrically from the entire
ensemble of available data, across all subjects or experiments.

\begin{appendix}
\section*{Appendix}\label{app}
\subsection{Additional details and
discussion for preliminary formulas}\label{secA1}

Formula (\ref{rep}) is an extension of the covariance kernel
representation in terms of eigenfunctions, given by
$\cov(X(t),X(s))=\sum\lambda_k \phi_k(t) \phi_k(s)$ [\citet{ash75}],
which itself is a stochastic process version of the classical
multivariate representation of a covariance matrix $C$ in terms of
its eigenvectors $e_k$ and eigenvalues $\lambda_k$, $C=\sum
\lambda_k e_k e_k'$. Specifically, using representation
(\ref{kl3}), one finds
$\cov(X^{(\nu_1)}(t),X^{(\nu_2)}(s))=\sum_{k,l} \cov(\xi_k,\xi_l)
\phi_k^{(\nu_1)}(t) \phi_k^{(\nu_2)}(s)$, and (\ref{rep}) follows
upon observing that $\cov(\xi_k,\xi_l)=\lambda_k$ for $k=l$ and
$=0$ for $k \ne l$.

Regarding the ``population differential equation''
$E\{X^{(1)}(t)-\mu^{(1)}(t)|\break X(t)\}=\beta(t)\{X(t)-\mu(t)\}$,
observe that
for any jointly normal random vectors $(U_1,U_2)$ with mean $0$ and
covariance matrix $C$ with elements $c_{11}, c_{12}, c_{21}=c_{12},
c_{22}$, it holds that $E(U_2|U_1)=(c_{21}/c_{11})U_1$. Applying
this to the jointly normal random vectors in (\ref{jt}) then
implies this population equation. The specific form for the function
$\beta$ in (\ref{beta}) is obtained by plugging in the specific
terms of the covariance matrix given on the right-hand side of
(\ref{jt}).

Applying (\ref{rep}), observing $\var(X(t))=\sum_k \lambda_k
\phi_k(t)^2$, and then taking the log-derivative leads to
$\frac{d}{dt}\log(\var(X(t)))=2 [\sum_k \lambda_k \phi_k(t)
\phi_k^{(1)}(t)]/[\sum_k \lambda_k \phi_k^2(t)]$, establishing the
last equality in representation (\ref{beta}).\vadjust{\goodbreak}

\subsection{Estimation procedures}\label{secA2}

Turning to estimation, in a first step we aggregate measurements
across all subjects into one ``big'' scatterplot and apply a
smoothing method that allows us to obtain the $\nu$th derivative of a
regression function from scatterplot data. For example, in the case
of local polynomial fitting, given a univariate density function
$\kappa_1$ and bandwidth $h_{\mu,\nu}$, one would minimize
%
%e21 ###
%
\begin{equation} \label{smooth1}
\sum_{i=1}^{n}\sum_{j=1}^{N_i}\kappa_1 \biggl(\frac{T_{ij}-t}
{h_{\mu,\nu}} \biggr) \Biggl\{Y_{ij}-\sum_{m=0}^{\nu+1}
\alpha_m(T_{ij}-t)^m \Biggr\}^2
\end{equation}
for each $t$ with respect to $\alpha_m$ for $m = 0,\ldots,\nu+1$,
from which one obtains $\hat{\mu}^{(\nu)}(t) =
\hat{\alpha}_{\nu}(t)\nu!$ [\citet{fan96}].

According to (\ref{kl4}), we will also need estimates of
$\frac{\partial^{\nu}}{\partial t^{\nu}}G(t,s)={G}^{(\nu,0)}$.
There are various techniques available for this task. Following
\citet{mull091}, to which we refer for further details, using again
local polynomial fitting, we minimize the pooled scatterplot of
pairwise raw covariances
%
%e22 ###
%
\begin{eqnarray}\label{smooth2}\quad
&&\sum_{i=1}^{n}\sum_{1 \leq j \neq l \leq{N_i}} \kappa_2
\biggl(\frac{T_{ij}-t}{h_{G,{\nu}}},\frac{T_{il}-s}{h_{G,{\nu
}}} \biggr)\nonumber\\[-8pt]\\[-8pt]
&&\qquad\hspace*{30.8pt}{}\times
\Biggl\{G_i(T_{ij},T_{il})-\Biggl(\sum_{m = 0}^{\nu+1}\alpha_{1m}
(T_{ij}-t)^m+\alpha_{21}(T_{il}-s)\Biggr) \Biggr\}^2\nonumber
\end{eqnarray}
for fixed $(t,s)$ with respect to $\alpha_{1m}$ and $\alpha_{21}$
for $m = 1,\ldots,\nu+1$, where $G_i(T_{ij}, T_{il}) =
(Y_{ij}-\hat{\mu}(T_{ij}))(Y_{il}-\hat{\mu}(T_{il}))$, $j \ne l$,
$\kappa_2$ is a kernel chosen as a bivariate density function, and
$h_{G,{\nu}}$ is a bandwidth. This leads to $\hat{G}^{(\nu,0)}(t,s)
= \hat{\alpha}_{1\nu}(t,s)\nu!$.

The pooling that takes place in the scatterplots for estimating the
derivatives of $\mu$ and of $G$ is the means to accomplish the
borrowing of information across the sample, which is essential to
overcome the limitations of the sparse sampling designs. We note
that the case of no derivative $\nu=0$ is always included, and
solving the eigenequations on the left-hand side of (\ref{kl4})
numerically\vspace*{1pt} for that case leads to the required estimates
$\hat{\lambda}_1, \hat{\lambda}_2,\ldots$ of the eigenvalues and
$\hat{\phi}_1, \hat{\phi}_2,\ldots$ of the eigenfunctions. The
estimates $\hat{\phi}_1^{(1)}, \hat{\phi}_2^{(1)},\ldots$ of the
eigenfunction derivatives are then obtained from the right-hand side
of (\ref{kl4}), plugging in the estimates for eigenfunctions and
eigenvalues, followed by a numerical integration step.

The plug-in estimates, $\hat{\beta}_K, \hat{G}_{z,K}, \hat{V}_K,
\hat{R}^2_K$, are then obtained from the corresponding
representations, (\ref{beta}), (\ref{Gz}), (\ref{v}), (\ref{R2}), by
including $K$ leading components in the respective sums. While for
theoretical analysis and asymptotic consistency one requires $K=K(n)
\rightarrow\infty$, the number of included eigen-terms $K$ in
practical data analysis can be chosen by various criteria, for example,
AIC/BIC based on marginal/conditional pseudo-likelihood or
thresholding of\vadjust{\goodbreak} the total variation explained by the included
components [\citet{mull091}]. One key feature of the covariance surface
smoothing step in~(\ref{smooth2}) is the exclusion of the diagonal
elements (for which $j=l$); the expected value for these elements
includes the measurement error variance $\sigma^2$ in addition to
the variance of the process. The difference between a smoother that
uses the diagonal elements only and the resulting diagonal from the
smoothing step (\ref{smooth2}) when no derivatives are involved can
then be used to find consistent estimates for the error variance
$\sigma^2$ [\citet{mull054}].

To obtain estimates for the derivatives of the trajectories $X_i$, a
realistic target is the conditional expectation
$E\{X^{(\nu)}_i(t)|Y_{i1},\ldots,Y_{iN_i}\}$. It turns out that this
conditional expectation can be consistently estimated in the case of
Gaussian processes by applying principal analysis by conditional
expectation (PACE) [\citet{mull054}]. For $\mathbf{X}_i =
(X_i(T_{i1}),\ldots, X_i(T_{iN_i}))^T$, $\mathbf{Y}_i =
(Y_{i1},\ldots,
Y_{iN_i})^T$, $\bolds{\mu}_i = (\mu(T_{i1}),\ldots, \mu(T_{iN_i}))^T$,
$\bolds\phi_{ik} =(\phi_k(T_{i1}), \ldots,\break\phi_k(T_{iN_i}))^T$, if
$\xi_{ik}$ and $\varepsilon_{ij}$ in (\ref{kl2}) are jointly
Gaussian, then by standard properties of the Gaussian distribution,
%
%e23 ###
%
\begin{equation}\label{score}
E(\xi_{ik} |\mathbf{Y}_i) = \lambda_k
\bolds\phi_{ik}^T\Sigma_{Y_i}^{-1}(\mathbf{Y}_i-\bolds\mu_i),
\end{equation}
where $\Sigma_{Y_i} = \cov(\mathbf{Y}_i,
\mathbf{Y}_i) = \cov(\mathbf{X}_i, \mathbf{X}_i)+\sigma^2
{I}_{N_i}$. This
implies $E(X^{(\nu)}_i(t)|Y_{i1},\break\ldots, Y_{iN_i})=\sum_{k=1}^\infty
E(\xi_{ik} |\mathbf{Y}_i) \phi_k^{(\nu)}(t)=\{\sum
_{k=1}^\infty
\lambda_k \phi_k^{(\nu)}(t)
\bolds\phi_{ik}^T\} \Sigma_{Y_i}^{-1}(\mathbf{Y}_i-\bolds\mu_i),
\nu=0,1$.
The unknown quantities can be estimated by simply plugging in the
variance, eigenvalue, eigenfunction and eigenfunction derivative
estimates discussed above, again coupled with truncating the number
of included components at $K$.

\subsection{Additional assumptions and auxiliary results}\label{secA3}

Denote the densities of $T$ and $(T_1, T_2)$ by $f_1(t)$, $f_2(t,
s)$, and define an interior domain by $\T=[a, b]$ with
$\T_\delta=[a-\delta, b+\delta]$ for some $\delta>0$. Regularity
conditions for the densities and the targeted moment functions as
well as their derivatives are as follows, where $\ell_1, \ell_2$ are
nonnegative integers:
\begin{enumerate}[(B2)]
\item[(B1)] $f_1^{(5)}(t)$ exists and is continuous on $\T_\delta$
with $f(t) > 0$,
$\frac{\partial^5}{\partial t^{\ell_1} \,\partial s^{\ell_2}}f_2(t,s)$
exists and is continuous on $\T_\delta^2$ for $\ell_1+\ell_2=5$;
\item[(B2)] $\mu^{(5)}(t)$ exists and is continuous on $\T_\delta$,
$\frac{\partial^5}{\partial t^{\ell_1}\, \partial s^{\ell_2}}G(t, s)$ exists
and is continuous on $\T_\delta^2$ for $\ell_1+\ell_2=5$.
\end{enumerate}

We say that a bivariate kernel function $\kappa_2$ is of order
$(\bolds\nu, \ell)$, where $\bolds\nu$ is a multi-index $\bolds
\nu=(\nu_1, \nu_2)$, if
%
%e24 ###
%
\begin{eqnarray}\label{order}
&&\int\int u^{\ell_1}
v^{\ell_2} K_2(u,v)\,du\,dv\nonumber\\[-8pt]\\[-8pt]
&&\qquad=
\cases{0, &\quad $0 \leq\ell_1+\ell_2 <\ell, \ell_1\neq\nu
_1, \ell_2 \neq\nu_2$, \cr
(-1)^{|{\bf\nu}|} \nu_1! \nu_2!, &\quad $\ell_1=\nu_1, \ell_2=\nu
_2$,\cr
\neq0, &\quad $\ell_1+\ell_2=\ell$,}\nonumber
\end{eqnarray}
where $|\bolds\nu|=\nu_1+\nu_2<\ell$. The univariate kernel
$\kappa_1$ is said to be of order $(\nu, \ell)$ for a univariate
$\nu=\nu_1$, if (\ref{order}) holds with $\ell_2=0$ on the
right-hand side, integrating only over the argument $u$ on the left-hand
side. For the kernel functions $\kappa_1$, $\kappa_2$ used in the
smoothing steps to obtain estimates for $\mu(t)$ and $\mu^{(1)}(t)$
in (\ref{smooth1}) and for $G(t,s)$ and $G^{(1,0)}(t,s)$ in
(\ref{smooth2}) we assume
\begin{enumerate}[(B3)]
\item[(B3)] Kernel functions $\kappa_1$ and $\kappa_2$ are
nonnegative with compact
supports, bounded and of order $(0, 2)$ and $((0, 0), 2)$,
respectively.
\end{enumerate}

The following lemma provides the weak $L^2$ convergence rate for
univariate and bivariate weighted averages defined below. For
arbitrary real functions $\theta\dvtx\Re^2\rightarrow\Re$ and
$\theta^\ast\dvtx\Re^4\rightarrow\Re$, define
$\tilde{\theta}(t)=E\{\theta(t_{ij}, Y_{ij})|T_{ij}=t\}$ and
$\tilde{\theta}^\ast(t)=E\{\theta^\ast(t_{ij}, t_{il}, Y_{ij},
Y_{il})|T_{ij}=t, T_{il}=s\}$, let
$\theta_\nu(t)=\tilde{\theta}^{(\nu)}(t) f_1(t)$ for a single index
$\nu$ and $\theta_{\bolds\nu}^\ast(t)=f_2(t, s)\frac{\partial
^{|\bolds\nu|}}{\partial^{\nu_1}\,\partial^{\nu_2}}\tilde{\theta
}^\ast(t, s)$
for a multi-index $\bolds\nu=(\nu_1, \nu_2)$, and define the
weighted kernel averages, employing bandwidths $h_1, h_2$,
%
%e26 ###
%e25 ###
%
\begin{eqnarray}
\label{kwa} \hat{\theta}(t)&=&\frac{1}{E(N)n h_1^{\nu+1}}\sum_{i=1}^n
\sum_{j=1}^{N_i} \theta(T_{ij}, Y_{ij})
\kappa_1 \biggl(\frac{t-T_{ij}}{h_1} \biggr),\\[-2pt]
\label{kwa1}
\hat{\theta}^\ast(t, s)&=&\frac{1}{E\{N(N-1)\} n
h_2^{|\bolds\nu|+2}}\sum_{i=1}^n \sum_{1\leq j\neq l \leq N_i}
\theta^\ast(T_{ij}, T_{il}, Y_{ij}, Y_{il})
\kappa_2\nonumber\\[-9pt]\\[-9pt]
&&\hspace*{152.4pt}{}\times\biggl(\frac{t-T_{ij}}{h_2}, \frac{s-T_{il}}{h_2}
\biggr).\nonumber
\end{eqnarray}
For establishing convergence results for the general weighted
averages (\ref{kwa}), assume that:
\begin{enumerate}[(B3$^\dag$)]
\item[(B2$^\dag$)] Derivatives $\tilde{\theta}^{(\ell)}(t)$ exist and
are continuous on $\T_\delta$, $\frac{\partial^5}{\partial
t^{\ell_1}\,
\partial s^{\ell_2}}\tilde{\theta}^\ast(t, s)$ exists and is
continuous on $\T_\delta^2$ for $\ell_1+\ell_2=\ell$.
\item[(B3$^\dag$)] The univariate kernel $\kappa_1$ is of order
$(\nu, \ell)$ and the bivariate kernel $\kappa_2$ is of order
$(\bolds\nu, \ell)$.
\end{enumerate}
\begin{lem} \label{lem1}
Under \textup{(A1)}, \textup{(B1)},
\textup{(B2$^\dag$)}, \textup{(B3$^\dag$)} and if $\max\{h_1, h_2\}
\rightarrow0$, $nh_1^{2\nu+1}\rightarrow\infty$ and $nh_2^{2|\bolds
\nu|+2}\rightarrow\infty$, as $ n \rightarrow\infty$,
%
%e27 ###
%
\begin{eqnarray}
\label{lem1-eq}
\|\hat{\theta}-\theta_\nu\|&=&O_p\biggl(\frac{1}{\sqrt{nh_1^{2\nu
+1}}}+h_1^{\ell-\nu}\biggr),\nonumber\\[-9pt]\\[-9pt]
\|\hat{\theta}^\ast-\theta_{\bolds\nu}^\ast\|_s&=&O_p\biggl(\frac
{1}{\sqrt{n}h_2^{|\bolds\nu|+1}}+h_2^{\ell-|\bolds\nu|}\biggr).\nonumber
\end{eqnarray}
\end{lem}

\subsection{Technical proofs}\label{secA4}

\mbox{}

\begin{pf*}{Proof of Theorem~\ref{thm1}}
Since $X, X^{(1)}$ are
jointly Gaussian processes, it is clear that $Z$ is Gaussian.\vadjust{\goodbreak}
Formula (\ref{beta}) for $\beta(t)$ is obtained by observing that
for joint Gaussian r.v.s, $E\{X^{(1)}(t)-\mu^{(1)}(t)|
X(t)-\mu(t)\}=[\cov\{X^{(1)}(t),X(t)\}/\var\{X(t)\}]\{X(t)-\mu(t)\}
$. Then the
properties of the functional principal component scores lead
directly to
%
%e28 ###
%
\begin{equation}\label{cov}\hspace*{28pt}
\cov\bigl\{X^{(1)}(t),X(t)\bigr\}=\sum
_{k=1}^\infty
\lambda_k\phi_k^{(1)}(t) \phi_k(t),\qquad
\var\{X(t)\}=\sum_{k=1}^\infty\lambda_k\phi_k(t)^2,
\end{equation}
whence
$\beta(t)=E\{X^{(1)}(t)-\mu^{(1)}(t)| X(t)-\mu(t)\}$. This
implies $E\{Z(t)\}=0$. According
to (\ref{de}),
$\cov\{Z(t),X(t)\}=\cov\{X^{(1)}(t),X(t)\}-\beta(t)\var\{X(t)\} =0$,
for all $t \in\mt$, using (\ref{cov}) and (\ref{beta}). This
implies the independence of $Z(t), X(t)$, due to the Gaussianity.
Next observe $\cov\{Z(t), Z(s)\}=\cov\{X^{(1)}(t)-\break\beta(t)X(t),
X^{(1)}(s)-\beta(s)X(s)\}$, from which one obtains the result by
straightforward calculation.
\end{pf*}
\begin{pf*}{Proof of Lemma~\ref{lem1}}
Since $N_i\stackrel{\mathrm{i.i.d.}}{\sim}N$ and $N$ is a bounded and
integer-valued random variable. Denote the upper bound by $M$. To
handle the one-dimensional case in (\ref{kwa}), we observe
\begin{eqnarray*}
\hat{\theta}(t)&=&\sum_{j=1}^{M} \frac{1}{E(N)}\frac{1}{n
h_1^{\nu+1}}\sum_{i=1}^n \theta(T_{ij}, Y_{ij})
\kappa_1 \biggl(\frac{t-T_{ij}}{h_1} \biggr)\mathbf{1}(N_i\geq j)\\
&\equiv&\sum_{j=1}^M \frac{1}{E(N)}\hat{\theta}_{\nu j}(t),
\end{eqnarray*}
where
$\mathbf{1}(\cdot)$ is the indication function. Note that for each $j$,
$\hat{\theta}_{\nu j}$ is obtained from an i.i.d. sample. Slightly
modifying the proof of Theorem 2 in \citet{hall841} for a
kernel of
order $(\nu, \ell)$ provides the weak convergence rate
$\|\hat{\theta}_{\nu j}-\theta_\nu P(N\geq
j)\|=O_p\{(nh_1^{2\nu+1})^{-1/2}+h_1^{\ell-\nu}\}$. It is easy to
check that $\sum_{j=1}^M P(N\geq j)=E(N)$, as $N$ is a positive
integer-valued random variable. Therefore,
\[
\|\hat{\theta}-\theta_\nu\|\leq\sum_{j=1}^M \frac{P(N\geq
j)}{E(N)} \biggl\|\frac{\hat{\theta}_{\nu j}}{P(N\geq j)}-\theta_\nu
\biggr\|=O_p\biggl(\frac{1}{\sqrt{nh_1^{2\nu+1}}}+h_1^{\ell-\nu}\biggr).
\]
Analogously, for the two-dimensional case in (\ref{kwa1}), let
\[
\hat{\theta}_{\bolds\nu j}^\ast=\frac{1}{n h_2^{|\bolds\nu
|+2}}\sum_{i=1}^n \theta^\ast(T_{ij}, T_{il}, Y_{ij}, Y_{il})
\kappa_2 \biggl(\frac{t-T_{ij}}{h_2}, \frac{s-T_{il}}{h_2}
\biggr)\mathbf{1}\{N_i\geq\max(j, l)\},
\]
and then $\hat{\theta}_{\bolds\nu}^\ast=\sum_{1\leq j\neq l \leq
M}[E\{N(N-1)\}]^{-1}
\hat{\theta}_{\mathbf{j}}^\ast$. Similarly to the above, one has
$\|\hat{\theta}_{\bolds\nu j}^\ast-\theta_{\bolds\nu} P\{N\geq
\max(j, l)\}\|_s=O_p\{(nh_2^{2|\bolds\nu|+2})^{-1/2}+h_2^{\ell
-|\bolds\nu|}\}$. Again it is easy to verify that $E\{N(N-1)\}=\sum
_{1\leq
j\neq l \leq M} P\{N\geq\max(j, l)\}$. The triangle inequality for
the $L^2$ distance entails $\|\hat{\theta}^\ast-\theta_{\bolds\nu
}\|_s
=O_p\{(nh_2^{2|\bolds\nu|+2})^{-1/2}+h_2^{\ell-|\bolds\nu|}\}$.
\end{pf*}
\begin{pf*}{Proof of Theorem~\ref{thm2}}
Note that the estimators
$\hat{\mu}$, $\hat{\mu}^{(1)}$, $\hat{G}$ and $\hat{G}^{(1, 0)}$ all
can be written as functions of the general averages defined in
(\ref{kwa}), (\ref{kwa1}). Slightly modifying the proof of Theorem 1
in \citet{mull091}, with $\sup$ rates replaced by the $L^2$ rates
given in Lemma~\ref{lem1}, then leads to the optimal $L^2$ weak
convergence rates for $\hat{\mu}^{\nu}$ and $\hat{G}^{(\nu, 0)}$ in
(\ref{thm2-eq1}).

For the convergence rate of $\hat{\phi}^{(1)}$, Lemma 4.3 in
\citet{bosq00} implies that
%
%e29 ###
%
\begin{equation}\label{bosq}
|\hat{\lambda}_k-\lambda_k|\leq\|\hat{G}-G\|_s,\qquad
\|\hat{\phi}_k-\phi_k\|\leq2\sqrt{2}\delta_k^{-1}\|\hat{G}-G\|_s,
\end{equation}
where $\delta_k$ is defined in (\ref{spacing}) and $\hat{G}$ is
an arbitrary estimate (or perturbation) of $G$. Denote the linear
operators generated from the kernels $G^{(1, 0)}$ and
$\hat{G}^{(1,0)}$ by $\mathbf{G}^{(1, 0)}$, respectively, $\hat
\mathbf{G}^{(1,0)}$. Noting that $\delta_k\leq\lambda_k$, one finds
%
%e30 ###
%
\begin{eqnarray} \label{pf1}
\bigl\|\phi_k^{(1)}-\phi_k^{(1)}\bigr\|&\leq& \frac{1}{\hat{\lambda}_k}\bigl\|
\hat\mathbf{G}^{(1,0)}\hat{\phi}_k-\mathbf{G}^{(1,0)}\phi_k\bigr\|+ \bigl\|
\mathbf{G}^{(1,0)}\phi_k\bigr\|\cdot
\biggl|\frac{1}{\hat{\lambda}_k}-\frac{1}{\lambda_k}\biggr| \nonumber\\
&\stackrel{p}{\asymp}&\frac{1}{\lambda}_k\bigl\{\bigl\|\hat
{G}^{(1,0)}-G^{(1,0)}\bigr\|_s+
\|\hat{\phi}_k-\phi_k\|\bigr\}+\frac{|\hat{\lambda}_k-\lambda
_k|}{\lambda_k^2}\\
&\stackrel{p}{\asymp}&
\frac{1}{\lambda}_k\biggl\{\bigl\|\hat{G}^{(1,0)}-G^{(1,0)}\bigr\|_s+\frac
{1}{\delta_k}\|\hat{G}-G\|_s\biggr\},\nonumber
\end{eqnarray}
which implies (\ref{thm2-eq2}).
\end{pf*}
\begin{pf*}{Proof of Theorem~\ref{thm3}}
From (\ref{bosq}) it is easy to see that
$|\hat\lambda_k-\lambda_k|=o_p(\|\hat\phi_k-\phi_k\|)$, and from both
(\ref{bosq}) and (\ref{pf1}) that
$\|\hat\phi_k-\phi_k\|=o_p(\|\hat\phi_k^{(1)}-\phi_k^{(1)}\|)$
uniformly in $k$. One then finds that $\|\hat{G}^{(1,1)}-G^{(1,1)}\|_s$
is bounded in probability by
\begin{eqnarray*}
&&\sum_{k=1}^K \bigl\|\hat\lambda_k \hat\phi_k^{(1)}\hat\phi^{(1)}_k
-\lambda_k\phi_k^{(1)}\phi_k^{(1)}\bigr\|_s+ \Biggl\|\sum_{k=K+1}^\infty
\lambda_k\phi_k^{(1)} \phi_k^{(1)}\Biggr\|_s\\
&&\qquad\stackrel{p}{\asymp}
\sum_{k=1}^K \lambda_k \bigl\|\hat{\phi}_k^{(1)}-\phi_k^{(1)}\bigr\|+\|R_{K,
1}^\ast\|,
\end{eqnarray*}
which implies that
$\|\hat{G}^{(1,1)}-G^{(1,1)}\|_s=O_p(\alpha_n+\|R_{K, 1}^\ast\|)$,
where $\alpha_n$ is defined in (\ref{rate}) and the remainder terms
in (\ref{rem}). Similar arguments lead to
$\|\hat{G}^{(1,1)}-G^{(1,1)}\|_u=O_p(\alpha_n+\|R_{K, 1}\|)$, noting
$\|R_{K, 1}^\ast\|\leq\|R_{K, 1}\|$ due to the Cauchy--Schwarz
inequality.

Regarding $\|\hat{\beta}_K-\beta\|$, one has
%Note that $\|G^{(\nu, \nu)}\|_u^2<\infty$ in (A3)
%implies that $\|\sum_{k=K+1}^\infty\lambda_k
%as $K\rightarrow\infty$.
%
\begin{eqnarray*}
\|\hat{\beta}_K\!-\!\beta\|&\!\leq\!& \frac{1}{\inf_t
G(t, t)} \Biggl(
\sum_{k=1}^K \bigl\|\hat\lambda_k \hat\phi_k\hat\phi^{(1)}_k -
\lambda_k\phi_k\phi_k^{(1)}\bigr\|+ \Biggl\|\sum_{k=1}^K \lambda_k\phi_k
\phi_k^{(1)}\Biggr\| \Biggr)\\
&&{}\! + \!\biggl\|\frac{\sum_{k=1}^K\hat\lambda_k\hat
\phi_k\hat\phi^{(1)}_k}{{\sum_{k=1}^K
\hat\lambda_k\hat\phi_k^2}{\sum_{k=1}^\infty\lambda_k\phi_k^2}}
\biggr\| \Biggl(\sum_{k=1}^K
\|\hat\lambda_k\hat\phi_k^2-\lambda_k\phi_k^2\|+\Biggl\|\sum
_{k=K+1}^\infty
\lambda_k\phi_k^2\Biggr\| \Biggr)\\
&\stackrel{p}{\asymp}&\sum_{k=1}^K \lambda_k
\bigl\|\hat\phi^{(1)}_k-\phi_k^{(1)}\bigr\|+\|R_{K,0}\|+\sqrt{\|R_{K,0}\|
\cdot
\|R_{K,1}\|},
\end{eqnarray*}
applying the Cauchy--Schwarz inequality to
$\|\sum_{k=K+1}^\infty\lambda_k\phi_k \phi_k^{(1)}\|$. Observing
$\sqrt{\|R_{K,0}\|\cdot\|R_{K,1}\|}\leq
(\|R_{K,0}\|+\|R_{K,1}\|\}/2$ yields
$\|\hat{\beta}_K-\beta\|=O_p(\alpha_n+\|R_{K,0}\|+\|R_{K,1}\|)$.

To study $\|\hat{G}_{z,K}-G_z\|_s$, we investigate the $L^2$
convergence rates of
\[
I_1=\Biggl\|\hat{\beta}_K \sum_{k=1}^K
\hat\lambda_k \hat\phi_k \hat\phi_k^{(1)}-\beta\sum
_{k=1}^\infty
\lambda_k \phi_k \phi_k^{(1)}\Biggr\|_s,\qquad I_2=\|\hat{\beta}_K
\hat{G}_K \hat\beta_K-\beta G\beta\|_s,
\]
where $\beta$ (resp.,
$\hat\beta_K$) and $\phi_k$ (resp., $\hat{\phi}_k$) share the same
argument, and we define $\hat{G}_K(t, s)=\sum_{k=1}^K
\hat\lambda\hat\phi_k(t) \hat\phi_k(s)$. In analogy to the above
arguments, $I_1\stackrel{p}{\asymp}
\|\hat{\beta}_K-\beta\|+\sum_{k=1}^K
\lambda_k\|\hat\phi_k^{(1)}-\phi_k^{(1)}\|+\sqrt{\|R_{K,0}\|\cdot
\|R_{K,1}\|}$, $I_2=\|\hat{\beta}_K-\beta\|+\sum_{k=1}^K
\lambda_k\|\hat\phi_k-\phi_k\|+\|R_{K,0}\|$. This leads to
$\|\hat{G}_{z,K}-G_z\|_s=O_p(\alpha_n+\|R_{K,0}\|+\|R_{K,1}\|)$. The
same argument also applies to $\|\hat{G}_{z,K}-G_z\|_u$. Next we
study $\|\hat{R}^2_K-R^2\|$ and find that $\|\sum_{k=1}^K
\hat{\lambda}_k
\hat{\phi}_k^{(\nu_1)}\hat{\phi}_k^{(\nu_2)}-\sum_{k=1}^\infty
\lambda_k\phi^{(\nu_1)} \phi^{(\nu_2)}\|\stackrel{p}{\asymp}\break
\sum_{k=1}^K \lambda_k (\|\hat{\phi}_k^{(\nu_1)}-\phi_k^{(\nu
_1)}\|+
\|\hat{\phi}_k^{(\nu_2)}-\phi_k^{(\nu_2)}\|)+
\|R_{K,\nu_1}\|+\|R_{K,\nu_2}\|=O_p(\alpha_n+\|R_{K,\nu_1}\|+\|
R_{K,\nu_2}\|)$.
Analogous arguments apply to $\|\hat{V}_K-V\|$, completing the
proof.
\end{pf*}
\end{appendix}

\section*{Acknowledgments}

We are grateful to two referees for helpful comments that led to
an improved version of the paper.

%suskaldyti doi

%
\printaddresses

\end{document}